\hsize=13.50cm
\vsize=18cm
\parindent=12pt   \parskip=0pt
\pageno=1

% deux jeux d'offsets. Le premier concerne la \magnification=1000
% (si l'image n'est pas agrandie)

\hoffset=15mm    % offset horizontal en \magnification=1000
\voffset=1cm    % offset vertical en \magnification=1000

% Le deuxieme concerne  la \magnification=1200
% (lorsque l'image est agrandie de 20%)

\ifnum\mag=\magstep1
\hoffset=-2mm   % offset horizontal en \magnification=1200
\voffset=.8cm   % offset horizontal en \magnification=1200
\fi

% --------------------------- Reglages --------------------------
% ----------- auquels il ne vaut mieux pas toucher --------------

\pretolerance=500 \tolerance=1000  \brokenpenalty=5000

% -------------------- Debut des macros privees -------------------
\catcode`\@=11
% --------------------- Les fontes --------------------------------

\font\eightrm=cmr8         \font\eighti=cmmi8
\font\eightsy=cmsy8        \font\eightbf=cmbx8
\font\eighttt=cmtt8        \font\eightit=cmti8
\font\eightsl=cmsl8        \font\sixrm=cmr6
\font\sixi=cmmi6           \font\sixsy=cmsy6
\font\sixbf=cmbx6

% Fontes AMS

\font\tengoth=eufm10       \font\tenbboard=msbm10
\font\eightgoth=eufm10 at 8pt      \font\eightbboard=msbm10 at 8 pt
\font\sevengoth=eufm7      \font\sevenbboard=msbm7
\font\sixgoth=eufm7 at 6 pt        \font\fivegoth=eufm5

 \font\tencyr=wncyr10       
\font\eightcyr=wncyr10 at 8 pt      
\font\sevencyr=wncyr10 at 7 pt      
\font\sixcyr=wncyr10 at 6 pt

% Pour que les accents se placent correctement en mode math en corps 8 et 6

\skewchar\eighti='177 \skewchar\sixi='177
\skewchar\eightsy='60 \skewchar\sixsy='60

% Nouvelles familles pour les maths

\newfam\gothfam           \newfam\bboardfam
\newfam\cyrfam

\def\tenpoint{%
  \textfont0=\tenrm \scriptfont0=\sevenrm \scriptscriptfont0=\fiverm
  \def\rm{\fam\z@\tenrm}%
  \textfont1=\teni  \scriptfont1=\seveni  \scriptscriptfont1=\fivei
  \def\oldstyle{\fam\@ne\teni}\let\old=\oldstyle
  \textfont2=\tensy \scriptfont2=\sevensy \scriptscriptfont2=\fivesy
  \textfont\gothfam=\tengoth \scriptfont\gothfam=\sevengoth
  \scriptscriptfont\gothfam=\fivegoth
  \def\goth{\fam\gothfam\tengoth}%
  \textfont\bboardfam=\tenbboard \scriptfont\bboardfam=\sevenbboard
  \scriptscriptfont\bboardfam=\sevenbboard
  \def\bb{\fam\bboardfam\tenbboard}%
 \textfont\cyrfam=\tencyr \scriptfont\cyrfam=\sevencyr
  \scriptscriptfont\cyrfam=\sixcyr
  \def\cyr{\fam\cyrfam\tencyr}%
  \textfont\itfam=\tenit
  \def\it{\fam\itfam\tenit}%
  \textfont\slfam=\tensl
  \def\sl{\fam\slfam\tensl}%
  \textfont\bffam=\tenbf \scriptfont\bffam=\sevenbf
  \scriptscriptfont\bffam=\fivebf
  \def\bf{\fam\bffam\tenbf}%
  \textfont\ttfam=\tentt
  \def\tt{\fam\ttfam\tentt}%
  \abovedisplayskip=12pt plus 3pt minus 9pt
  \belowdisplayskip=\abovedisplayskip
  \abovedisplayshortskip=0pt plus 3pt
  \belowdisplayshortskip=4pt plus 3pt
  \smallskipamount=3pt plus 1pt minus 1pt
  \medskipamount=6pt plus 2pt minus 2pt
  \bigskipamount=12pt plus 4pt minus 4pt
  \normalbaselineskip=12pt
  \setbox\strutbox=\hbox{\vrule height8.5pt depth3.5pt width0pt}%
  \let\bigf@nt=\tenrm       \let\smallf@nt=\sevenrm
  \normalbaselines\rm}

\def\eightpoint{%
  \textfont0=\eightrm \scriptfont0=\sixrm \scriptscriptfont0=\fiverm
  \def\rm{\fam\z@\eightrm}%
  \textfont1=\eighti  \scriptfont1=\sixi  \scriptscriptfont1=\fivei
  \def\oldstyle{\fam\@ne\eighti}\let\old=\oldstyle
  \textfont2=\eightsy \scriptfont2=\sixsy \scriptscriptfont2=\fivesy
  \textfont\gothfam=\eightgoth \scriptfont\gothfam=\sixgoth
  \scriptscriptfont\gothfam=\fivegoth
  \def\goth{\fam\gothfam\eightgoth}%
  \textfont\cyrfam=\eightcyr \scriptfont\cyrfam=\sixcyr
  \scriptscriptfont\cyrfam=\sixcyr
  \def\cyr{\fam\cyrfam\eightcyr}%
  \textfont\bboardfam=\eightbboard \scriptfont\bboardfam=\sevenbboard
  \scriptscriptfont\bboardfam=\sevenbboard
  \def\bb{\fam\bboardfam}%
  \textfont\itfam=\eightit
  \def\it{\fam\itfam\eightit}%
  \textfont\slfam=\eightsl
  \def\sl{\fam\slfam\eightsl}%
  \textfont\bffam=\eightbf \scriptfont\bffam=\sixbf
  \scriptscriptfont\bffam=\fivebf
  \def\bf{\fam\bffam\eightbf}%
  \textfont\ttfam=\eighttt
  \def\tt{\fam\ttfam\eighttt}%
  \abovedisplayskip=9pt plus 3pt minus 9pt
  \belowdisplayskip=\abovedisplayskip
  \abovedisplayshortskip=0pt plus 3pt
  \belowdisplayshortskip=3pt plus 3pt
  \smallskipamount=2pt plus 1pt minus 1pt
  \medskipamount=4pt plus 2pt minus 1pt
  \bigskipamount=9pt plus 3pt minus 3pt
  \normalbaselineskip=9pt
  \setbox\strutbox=\hbox{\vrule height7pt depth2pt width0pt}%
  \let\bigf@nt=\eightrm     \let\smallf@nt=\sixrm
  \normalbaselines\rm}

\tenpoint

% Definition des petites capitales qui reagissent au \tenpoint et \eightpoint
% La syntaxe est celle d'un changement de fonte :
% {\pc FERMAT}, {\pc EUCLIDE} et {\pc G\"ODEL}.

\def\pc#1{\bigf@nt#1\smallf@nt}         \def\pd#1 {{\pc#1} }

%---------------------- dactylographie francaise -----------------

\catcode`\;=\active
\def;{\relax\ifhmode\ifdim\lastskip>\z@\unskip\fi
\kern\fontdimen2  -1.2 \fontdimen3 \string;}

\catcode`\:=\active
\def:{\relax\ifhmode\ifdim\lastskip>\z@\unskip\fi\penalty\@M\ \fi\string:}

\catcode`\!=\active
\def!{\relax\ifhmode\ifdim\lastskip>\z@
\unskip\fi\kern\fontdimen2  -1.1 \fontdimen3 \string!}

\catcode`\?=\active
\def?{\relax\ifhmode\ifdim\lastskip>\z@
\unskip\fi\kern\fontdimen2  -1.1 \fontdimen3 \string?}

\def\^#1{\if#1i{\accent"5E\i}\else{\accent"5E #1}\fi}
\def\"#1{\if#1i{\accent"7F\i}\else{\accent"7F #1}\fi}

\frenchspacing

% ---------------- Le format de sortie ----------------------------
% Haut et bas de page

\newtoks\auteurcourant      \auteurcourant={\hfil}
\newtoks\titrecourant       \titrecourant={\hfil}

\newtoks\hautpagetitre      \hautpagetitre={\hfil}
\newtoks\baspagetitre       \baspagetitre={\hfil}

\newtoks\hautpagegauche
\hautpagegauche={\eightpoint\rlap{\folio}\hfil\the\auteurcourant\hfil}
\newtoks\hautpagedroite
\hautpagedroite={\eightpoint\hfil\the\titrecourant\hfil\llap{\folio}}

\newtoks\baspagegauche      \baspagegauche={\hfil}
\newtoks\baspagedroite      \baspagedroite={\hfil}

\newif\ifpagetitre          \pagetitretrue

% \nopagenumbers : c'est un peu violent, mais ça marche. Alors ...

\headline={\ifpagetitre\the\hautpagetitre
\else\ifodd\pageno\the\hautpagedroite\else\the\hautpagegauche\fi\fi}

\footline={\ifpagetitre\the\baspagetitre\else
\ifodd\pageno\the\baspagedroite\else\the\baspagegauche\fi\fi
\global\pagetitrefalse}

% Redefinition de \raggedbottom pour avoir plus de mou en bas de page
% (necesssaire quand il y a beaucoup de grumeaux, des grosses
% formules centrees et pas beaucoup de texte entre)

\def\raggedbottom{\topskip 10pt plus 36pt\r@ggedbottomtrue}

% ------------------ Macros de mise en page ---------------------

% Un point-tiret

\def\pointir{\unskip . --- \ignorespaces}

% Macros Bigbreak et \Medbreak pour que les blancs verticaux ne s'ajoutent pas

\def\Bigbreak{\vskip-\lastskip\bigbreak}
\def\Medbreak{\vskip-\lastskip\medbreak}

% Texte centre dans une boite : le resultat est le plus petit rectangle
% qui contient le texte. Ce resultat est
% place dans une boite centree.
% La syntaxe est celle d'un tableau \`a une colonne

\def\ctexte#1\endctexte{%
  \hbox{$\vcenter{\halign{\hfill##\hfill\crcr#1\crcr}}$}}

% Titres centres (en gras)

\long\def\ctitre#1\endctitre{%
    \ifdim\lastskip<24pt\vskip-\lastskip\bigbreak\bigbreak\fi
  		\vbox{\parindent=0pt\leftskip=0pt plus 1fill
          \rightskip=\leftskip
          \parfillskip=0pt\bf#1\par}
    \bigskip\nobreak}

\long\def\section#1\endsection{%
\vskip 0pt plus 3\normalbaselineskip
\penalty-250
\vskip 0pt plus -3\normalbaselineskip
\Bigbreak
\message{[section \string: #1]}{\bf#1\unskip}\pointir}

\long\def\sectiona#1\endsection{%
\vskip 0pt plus 3\normalbaselineskip
\penalty-250
\vskip 0pt plus -3\normalbaselineskip
\Bigbreak
\message{[sectiona \string: #1]}%
{\bf#1}\medskip\nobreak}

\long\def\subsection#1\endsubsection{%
\Medbreak
{\it#1\unskip}\pointir}

\long\def\subsectiona#1\endsubsection{%
\Medbreak
{\it#1}\par\nobreak}

\def\rem#1\endrem{%
\Medbreak
{\it#1\unskip} : }

\def\remp#1\endrem{%
\Medbreak
{\pc #1\unskip}\pointir}

\def\rema#1\endrem{%
\Medbreak
{\it #1}\par\nobreak}

\def\newparwithcolon#1\endnewparwithcolon{
\Medbreak
{#1\unskip} : }

\def\newparwithpointir#1\endnewparwithpointir{
\Medbreak
{#1\unskip}\pointir}

\def\newpara#1\endnewpar{
\Medbreak
{#1\unskip}\smallskip\nobreak}

%---------------------------------------------
% enonces de theoremes avec numerotation apres
% #1 = THEOREME, COROLLAIRE, etc.
% #2 = numero (par exemple 3, 3.1, etc.)
% #3 = l'enonce du th proprememnt dit.

\long\def\th#1 #2\enonce#3\endth{%
   \Medbreak
   {\pc#1} {#2\unskip}\pointir{\it #3}\medskip}

\long\def\tha#1 #2\enonce#3\endth{%
   \Medbreak
   {\pc#1} {#2\unskip}\par\nobreak{\it #3}\medskip}

%---------------------------------------------
% enonces de theoremes avec numerotation d'abord
% #1 = numero (par exemple 3, 3.1, etc.)
% #2 = THEOREME, COROLLAIRE, etc.
% #3 = l'enonce du th proprememnt dit.

\long\def\Th#1 #2 #3\enonce#4\endth{%
   \Medbreak
   #1 {\pc#2} {#3\unskip}\pointir{\it #4}\medskip}

\long\def\Tha#1 #2 #3\enonce#4\endth{%
   \Medbreak
   #1 {\pc#2} #3\par\nobreak{\it #4}\medskip}
%---------------------------------------------

% les differents retraits, voir aussi \item

\def\decale#1{\smallbreak\hskip 28pt\llap{#1}\kern 5pt}
\def\decaledecale#1{\smallbreak\hskip 34pt\llap{#1}\kern 5pt}
\def\puce{\smallbreak\hskip 6pt{$\scriptstyle\bullet$}\kern 5pt}

% -----------------------------------------------------------------
% -------------- ce que Knuth n'a pas fait ------------------------
% -----------------------------------------------------------------

% Un \displaylines qui numerote à droite.
% La syntaxe est la meme que celle de \eqalignno

\def\displaylinesno#1{\displ@y\halign{
\hbox to\displaywidth{$\@lign\hfil\displaystyle##\hfil$}&
\llap{$##$}\crcr#1\crcr}}

% Un \displaylines qui numerote a gauche.
% La syntaxe est la meme que celle de \leqalignno

\def\ldisplaylinesno#1{\displ@y\halign{
\hbox to\displaywidth{$\@lign\hfil\displaystyle##\hfil$}&
\kern-\displaywidth\rlap{$##$}\tabskip\displaywidth\crcr#1\crcr}}

% Un \eqalign qui accepte plusieurs alignements verticaux
% Motif : \hfil ** & ** \hfil & \hfil ** & ** \hfill, etc.

\def\eqalign#1{\null\,\vcenter{\openup\jot\m@th\ialign{
\strut\hfil$\displaystyle{##}$&$\displaystyle{{}##}$\hfil
&&\quad\strut\hfil$\displaystyle{##}$&$\displaystyle{{}##}$\hfil
\crcr#1\crcr}}\,}

% Systeme d'equations precede d'une accolade.
% Copi\'e sur \eqalign, on s'en sert comme une matrice
% syntaxe : signe & coef & inconnue
% les coef sont justifi\'es \`a droite (\hfil coef)
% et les inconnues \`a gauche (inconnue\hfil)
% attention : un seul & ou deux && avant le signe =
% selon la justification choisie !
% Exemple : $$\system{
%             &2 &x &- &3 & y & = &&  -5 \cr
%            -&  &x &+ &  & y & = &&   6 \cr
%           }$$

\def\system#1{\left\{\null\,\vcenter{\openup1\jot\m@th
\ialign{\strut$##$&\hfil$##$&$##$\hfil&&
        \enskip$##$\enskip&\hfil$##$&$##$\hfil\crcr#1\crcr}}\right.}

% pour avoir des messages raisonnables avec les lettres accentu\'ees

\let\@ldmessage=\message

\def\message#1{{\def\pc{\string\pc\space}%
                \def\'{\string'}\def\`{\string`}%
                \def\^{\string^}\def\"{\string"}%
                \@ldmessage{#1}}}

% ----------------- Divers gadgets --------------------------------

% Pour se rendre la vie facile : \up{er}, \up{i\`eme}, n\up{0}, etc.

\def\up#1{\raise 1ex\hbox{\smallf@nt#1}}

% Utilisation : \cf. \etc.

\def\qed{\raise -2pt\hbox{\vrule\vbox to 10pt{\hrule width 4pt
                 \vfill\hrule}\vrule}}

\def\virg{\raise .4ex\hbox{,}}   % virgule après une fraction

 % point-virgule de ponctuation en maths

\def\build#1_#2^#3{\mathrel{
\mathop{\kern 0pt#1}\limits_{#2}^{#3}}}

% Entoure #2 d'un filet. Le filet est ecarte tout autour de #1
% Syntaxe \boxit{5pt}{...}. La ligne de base n'est pas perdue.

\def\boxit#1#2{%
\setbox1=\hbox{\kern#1{#2}\kern#1}%
\dimen1=\ht1 \advance\dimen1 by #1 \dimen2=\dp1 \advance\dimen2 by #1
\setbox1=\hbox{\vrule height\dimen1 depth\dimen2\box1\vrule}%
\setbox1=\vbox{\hrule\box1\hrule}%
\advance\dimen1 by .6pt \ht1=\dimen1
\advance\dimen2 by .6pt \dp1=\dimen2  \box1\relax}

% -----------------------------------------------------------------
% fin des macros privees
% -----------------------------------------------------------------

\catcode`\@=12

% pour qu'il la ferme
\showboxbreadth=-1  \showboxdepth=-1

 %%%%%
\input amssym.def
\input amssym.tex

\magnification=\magstep1
\hsize=17,5truecm
\vsize=25.5truecm
\hoffset=-0.9truecm
\voffset=-0.8truecm
%\nopagenumbers    %pagenumbers
\topskip=1truecm
\footline={\tenrm\hfil\folio\hfil}
\raggedbottom
\abovedisplayskip=3mm %Reduction of space between text and formulae
\belowdisplayskip=3mm
\abovedisplayshortskip=0mm
\belowdisplayshortskip=2mm
\normalbaselineskip=12pt  %This is default. Do NOT change it!!
\normalbaselines

%%%%%%%%%%%%%%%%%%%%%%

%---------------------------------------------------------

\def\diagram#1{\def\normalbaselines{\baselineskip=0pt\lineskip=5pt}
\matrix{#1}}

\def\hfl#1#2#3{\smash{\mathop{\hbox to#3{\rightarrowfill}}\limits
^{\scriptstyle#1}_{\scriptstyle#2}}}

\def\gfl#1#2#3{\smash{\mathop{\hbox to#3{\leftarrowfill}}\limits
^{\scriptstyle#1}_{\scriptstyle#2}}}

\def\vfl#1#2#3{\llap{$\scriptstyle #1$}
\left\downarrow\vbox to#3{}\right.\rlap{$\scriptstyle #2$}}

%------------------------------------------------------------------

\def\ra{\rightarrow}

\def\ov{\overline}
\def\Q{{\bf Q}}
\def\Z{{\bf Z}}
\def\P{{\bf P}}
\def\F{{\bf F}}
\def\G{{\bf G}}
\def\R{{\bf R}}
\def\Z{{\bf Z}}
\def\C{{\bf C}}
\def\lra{\longrightarrow}
\def\Br{{\rm Br}}
\def\A{{\bf A}}
\def\O{{\cal O}}
\def\X{{\cal X}}
\def\Y{{\cal Y}}
\def\S{{\cal S }}
\def\inv{{\rm inv}}
\def\Sha{{\cyr X}}
\def\x{{\bf x}}

% This is plain TEX

\centerline{\bf Beyond the Manin obstruction}

\bigskip
\centerline{\bf Alexei N. Skorobogatov}
\bigskip
{\bf 1. Introduction}
\medskip

Let $X$ be a smooth variety over a field
$k$, and $\Br(X)$ be the (cohomological) 
Brauer--Grothendieck group of $X$, 
$\Br(X)=H^2(X,\G_m)$.
Let $\overline k$ be a separable closure of $k$, 
$G=Gal(\overline k/k)$, and let
$\overline X=X\times _k \overline k$.
Let $\Br_0(X)$ (resp. $\Br_1(X)$) be the image of the
natural map $\Br(k)\ra\Br(X)$ (resp. the kernel of the natural map $\Br(X)\ra \Br(\overline X)^G$).
For any extension $K/k$ let
$$X(K)\times \Br(X) \ra \Br(K),\ \  (A,P)\mapsto A(P)$$
be the pairing obtained by specializing elements of $\Br(X)$
at $K$-rational points. Now
let $k$ be a number field, $\Omega$ the set of places
of $k$, $\A_k$ the ad\`ele ring of $k$. Let
$$\inv_v: \Br(k_v)\hookrightarrow \Q/\Z$$
be the local invariant. In [M] Manin considered
the pairing between $\Br(X)$ and $X(\A_k)$ with values in $\Q/\Z $ 
given by
$$\sum_{v\in \Omega}\inv_v(A(P_v)),$$
for $A\in \Br(X)$, $\{P_v\}\in X(\A_k)$, the sum which 
is well known to
be finite (for almost all $v$ the specialization
$A(P_v)$ belongs to $\Br(\O_v)=0$).  
By the global reciprocity this pairing 
is trivial on the algebras coming from $\Br(k)$, so it
could be regarded as a pairing
$$\Br(X)/\Br_0(X)\times  X(\A_k)\ra \Q/\Z.$$
We shall call this {\it the Brauer--Manin pairing}.
Let us define
$X(\A_k)^\Br$ as ``the right kernel" of this pairing,
that is, the subset of points of $X(\A_k)$ orthogonal to 
all elements of $\Br(X)$. Manin made an important
observation that by the global reciprocity law the image of
$X(k)$ under the diagonal embedding $X(k)\hookrightarrow
X(\A_k)$ is contained in $X(\A_k)^\Br$. A variety
$X$ such that $X(\A_k)\not=\emptyset$ whereas $X(k)=\emptyset$
is {\it a counterexample to the Hasse principle}. Such
a counterexample is accounted for by {\it the Manin obstruction}
if already $X(\A_k)^\Br$ is empty. For a long time
most known counterexamples to the Hasse principle 
could be explained by means of the Manin obstruction 
(to the best of my knowledge the case of the
Bremner--Lewis--Morton curve $3x^4+4y^4-19z^4=0$
remains undecided).
Recently
Sarnak and L.Wang [SW] showed that the Manin obstruction
is not the only obstruction to the Hasse principle for
smooth hypersurfaces of degree 1130 
in $\P^4_\Q$ if one assumes Lang's
conjecture that $X(\Q)$ is finite if $X_\C$ is hyperbolic.

The aim of this note is to
construct a smooth proper surface over $k=\Q$ 
of Kodaira dimension 0 which is
a counterexample to the Hasse principle but for which
the Manin obstruction is not sufficient to explain the
absence of $\Q$-rational points. We exploit the same kind
of surfaces
which has been recently used to produce counterexamples
to a conjecture of Mazur [AA]. These surfaces are quotients 
of a product of two curves of genus one by a 
fixed point free involution. 

Although in our example the Manin obstruction fails 
to provide a finite decision process for determining
the existence of rational points there still is such 
a process. We propose a refinement of the Manin
obstruction and show that for surfaces of our type it 
is the only obstruction to the Hasse principle. 
To define it we use a combination of the Manin
obstruction and a descent very similar to the classical
descent on elliptic curves. The point is that unlike
in that classical case in the case of surfaces
the Brauer group can become substantially bigger
after passing to a finite unramified covering, thus
the Manin obstruction can become finer. The refined
obstruction depends on the choice of a finitely
generated submodule of $Pic(\ov X)$, and it can 
give something non-trivial
only if there is such a ``non-trivial" submodule. This
obstruction is thus unable to explain the example of 
Sarnak and L.Wang where $Pic(\ov X)=\Z$.

We construct our example in Section 2. A refinement
of the Manin obstruction is defined in Section 3. There
we also formulate the
underlying descent statement, which is then 
proved in Section 4. 
The construction of the example 
relies on the existence of elliptic
curves over $\Q$ with no rational 2-torsion and an
element of exact order 4 in the Tate--Shafarevich group.
In the appendix to this paper S. Siksek writes down
explicitly an everywhere locally soluble 4-covering 
of the curve $y^2=x^3-1221$, and shows that 
it does represent an element of exact order 4.

\bigskip
{\bf 2. A counterexample to the Hasse principle 
not accounted for by the Manin obstruction}

\bigskip
Let $J$ be an elliptic curve over ${\bf Q}$ whose Tate-Shafarevich group ${\cyr X}(J)$ 
contains an element of order 4, and such that 
$J$ contains no non-trivial point of order 2 defined 
over $\Q$.
For example, take $J$ to be the curve
$$y^2=x^3-1221.$$
Then $J(\Q)=0$ and ${\cyr X}(J)=\Z/4\times \Z/4$ (the
last property is conditional on 
the Birch - Swinnerton-Dyer conjecture, [GPZ],
Tables 3 and 4). We shall only use the 
(unconditional) result proved in
the appendix to this paper 
that ${\cyr X}(J)$ contains an element of exact order 4.
Let $C'$ be a principal homogeneous space 
under $J$ whose class $[C']\in H^1(\Q,J)$ belongs to 
${\cyr X}(J)$ and has exact order 4.

Let $C$ be a principal homogeneous space under $J$ 
such that $[C]=2[C']$. Let $\xi:C'\ra C$ be the corresponding
unramified Galois covering with Galois group $J[2]$.
In particular, the inverse image on the Jacobians,
$\xi^*:J\ra J$, is multiplication by 2.

By a lemma of Swinnerton-Dyer 
([B/SwD], Lemmas 1 and 2, [C], IV, Thm. 1.3) 
any element of order 2 in ${\cyr X}(J)$
can be represented as a double cover
of ${\bf P}^1_{\bf Q}$. In particular, $C$ can be given by
the equation
$$u^2=g(t)$$
for some polynomial $g(t)$ of degree 4 with integral 
coefficients (see (A.2) for an
explicit expression of $g(t)$). 
Let $\sigma:C\ra C$ be the corresponding
hyperelliptic involution. In particular, $\sigma_*$ 
acts on the
Jacobian $J$ of $C$ as multiplication by $-1$.

Let us fix once and for all two monic quadratic polymomials
$p(x)$ and $q(x)$ with integer coefficients such that
$Res(p(x),q(x))=\pm 1$, and such that both $p(x)$ and $q(x)$
take only positive values on $\Q$.
Let $D\subset {\bf P}^3_{\bf Q}$ be the smooth proper
curve of genus one given by its affine equations
$$y^2=p(x),\ z^2=q(x).$$
This curve $D$ has obvious rational points at infinity.
We fix one of them, say, $Q$, and take it 
as the neutral element of the group law on $D$. 

Let $\rho:D\ra D$
be the involution sending $(x,y,z)$ to $(x,-y,-z)$. 
One sees
easily that $\rho$ has no fixed point. Then $D'=D/\rho$ is
an elliptic curve given by
$$w^2=p(x)q(x).$$
Let $\psi:D\ra D'$ be the natural surjection. Take
$Q'=\psi(Q)$ as the neutral element for the group law 
on $D'$. Then $\psi$ becomes an isogeny whose kernel is 
generated by a point of order 2. 
In particular, the action of $\rho_*$
on the Jacobian of $D$ is trivial. (After the choice of
$Q$ and $Q'$ we can identify $D$ and $D'$ with their 
Jacobians.)

Let $X$ be the quotient of $Y=C\times D$ by the fixed point 
free involution $(\sigma,\rho)$, $f:Y\ra X$. This is a
surface classically known as hyperelliptic (or bielliptic), see [Sh], Ch. VII.8, [B], Ch. VI. Its geometric invariants are $\kappa=0$, $p_g=0$, $q=1$, $(K_X^2)=0$, $b_1=b_2=2$. 
(Recall that surfaces with such invariants
together with $K3$, Enriques and Abelian surfaces
exhaust all surfaces of Kodaira dimension zero.)
An affine
model of $X$ can be given by equations
$$y^2=g(t)p(x),\ z^2=g(t)q(x).\leqno{(2.1)}$$
Finally, let $Y'=C'\times D$, and let $\pi:Y'\ra D$ be the 
second projection.
Let $f':Y'\ra X$ be the composition of the unramified
covering $(\xi, Id): Y'\ra Y$ with $f$.

\medskip

{\bf Theorem.} (a) {\it We have} $X({\bf Q})=\emptyset$. 

(b) {\it $f'^{*}(\Br(X))\subset \Br(Y')$ is contained in 
$\pi^*(\Br(D))$. Therefore for any $R\in D({\Q})$, and any
$\{P_v\}\in C'(\A_{\Q})$, the map $f'$ sends 
$\{(P_v,R)\}\in Y'({\A}_{\Q})$ to
$X({\A}_{\Q})^{\Br}$, in particular we have}
$X({\A}_{\Q})^{\Br}\not=\emptyset$.

\medskip

{\it Proof.} (a) An easy valuation argument 
(the same as in [AA], Prop. 5.1) shows that
for any prime $p$ and 
any $\Q_p$-rational point of $X$ for which $yz\not=0$ 
the $p$-adic valuation
$val_p(g(t))$ is even. Thus $g(t)=\pm 1$
modulo squares in $\Q^*$. This extends to the whole of $X$
(see {\it loc.cit.}), and we get a decomposition
$$X({\bf Q})=f(Y({\bf Q}))\cup 
f^-(Y^-({\bf Q})),$$ 
where $f^-:Y^-=C^-\times D^-\ra X$ is a ``twisted form" of $f$, and $C^-$ (resp. $D^-$) is obtained by inverting
the sign of $g(t)$ (resp. of $p(x)$ and $q(x)$). Note that
$Y=C\times D$ clearly has no
rational point since $C({\bf Q})$ is empty. The curve
$D^-$ given by the equations
$y^2+p(x)=z^2+q(x)=0$ has no real point
by the positivity condition on $p(x)$ and $q(x)$, hence 
$Y^-$ has no real point.
This completes the proof of (a).

\medskip 

(b) The second statement of (b) follows from the first
one:
by projection formula it is enough to show that 
$\{(P_v,R)\}$ is Brauer--Manin orthogonal to 
$f'^{*}(\Br(X))$, and this follows from the first statement
by the global reciprocity law. 

Our first goal will be to study $f^*(\Br(X))$. 

Consider the Hochschild--Serre
spectral sequence ([Mi2], III.2.20):
$$H^p(G,H^q(\overline X,\G_m))\Rightarrow H^{p+q}(X,\G_m)
\leqno{(2.2)}$$
When $H^0(\ov X,\G_m)=\ov k^*$, for example when $X$ is proper, the
exact sequence of low degree terms writes down as follows:
$$0\ra \Br_0(X)\ra \Br_1(X)\ra H^1(G,Pic(\overline X))\ra
H^3(\Q,\G_m)=0,$$
where the last equality is provided by the class field theory 
(and holds also with $\Q$ replaced by any number
field or any local field).
We shall employ the notation 
$r:\Br_1(X)\ra H^1(G,Pic(\overline X))$ for the 
corresponding canonical map.

For surfaces the structure of $\Br(\ov X)$ was
determined by Grothendieck ([G], II, Cor. 3.4, III, (8.12)):
the divisible subgroup $\Br(\ov X)_{div}$ is
isomorphic to $(\Q/\Z)^{b_2-\rho}$ as an abelian group,
and the quotient $\Br(\ov X)/\Br(\ov X)_{div}$
is isomorphic to $Hom(NS(\overline X)_{tors},\Q/\Z)$
as a $G$-module. Here $b_2=rk(H^2(X_{\C},\Z))$, 
$\rho=cork(Pic(\overline X)\otimes \Q/\Z)$,
and $NS(\overline X)$ is the N\'eron--Severi group of
$\overline X$, finitely generated and isomorphic to
the quotient of $Pic(\overline X)$ by its divisible 
subgroup. Since in our case the geometric genus of $X$ is
zero,
we have $b_2=\rho$, hence $\Br(\overline X)$ is dual
to the torsion subgroup of $NS(\overline X)$.

We now analyse the structure of the $G$-module 
$Pic(\overline X)$ and of the map $f^*:Pic(\overline X)\ra
Pic(\ov Y)$. Consider the 
following commutative diagram:
$$\diagram{
C&\buildrel{\pi_1}\over{\hbox to 8mm{\leftarrowfill}} 
&Y&\buildrel{\pi_2}\over{\hbox to 8mm{\rightarrowfill}} 
&D\cr
\vfl{}{}{3mm}&& \vfl{f}{}{3mm}&& \vfl{\psi}{}{3mm}\cr
\P^1_\Q&\buildrel{}\over{\hbox to 8mm{\leftarrowfill}}   &X&\buildrel{}\over{\hbox to 8mm{\rightarrowfill}} 
&D'\cr
}$$
Here $X\ra D'$ is the Albanese map ([AA], Prop. 3.1 (3)). 
Let $\psi^*:D'\ra D$ be the isogeny
dual to $\psi: D\ra D'$. 
By the choice of $Q'\in D'(\Q)$ the $G$-module 
$D'(\ov k)$ is identified with $Pic^0(\ov X)$. 
Let $\Gamma$ be the cyclic group of automorphisms of 
$Y$ generated by
${(\sigma, \rho)}$. Since $f^*: Pic(\ov X)
\ra Pic(\overline Y)$ factors through the inclusion 
$Pic(\overline Y)^{\Gamma}\hookrightarrow
Pic(\overline Y)$ we get
a commutative diagram of $G$-modules with exact rows:
$$\diagram{
0 & \ra & D'(\ov k) & \ra & Pic(\overline X)&\ra 
& NS(\overline X)&\ra & 0\cr
& & \vfl{(0,\psi^*)}{}{3mm}& & \vfl{f^*}{}{3mm}& & \vfl{f^*}{}{3mm}\cr
0 & \ra & J[2](\ov k)\times D(\ov k) & \ra &
Pic(\overline Y)^{\Gamma}&\ra
& NS(\overline Y)^{\Gamma}\cr
& & \vfl{}{}{3mm}& & \vfl{}{}{3mm}& & \vfl{}{}{3mm}\cr
0 & \ra & J(\ov k)\times D(\ov k) & \ra & 
Pic(\overline Y)&\ra
& NS(\overline Y)&\ra&0\cr
}\leqno{(2.3)}$$
We used the fact that $(\sigma_*, \rho_*)$ acts on 
the product of two Jacobians
$J(\ov k)\times D(\ov k)$ as multiplication by $(-1,1)$. 
\bigskip
{\bf Lemma 1.} {\it There is an isomorphism of $G$-modules}
$NS(\ov X)_{tors}=J[2](\ov k)$.
\medskip
{\it Proof.} Let us denote the kernel of $f^*:NS(\ov X)\ra NS(\ov Y)$ by $K$.
Since $NS(\ov Y)$ is torsion-free, we have 
$NS(\ov X)_{tors}= K_{tors}$. Thus it will be 
enough to show that the $G$-modules $K$ and $J[2]$
are isomorphic.
The Hochschild--Serre spectral sequence
$$H^p(\Gamma, H^q(\overline Y, \G_m))
\Rightarrow H^{p+q}(\ov X,\G_m)$$
yields an exact sequence
$$0\ra H^1(\Gamma,\ov \Q^*)\ra Pic(\ov X) 
\buildrel{f^*}\over{\hbox to 6mm{\rightarrowfill}} 
Pic(\ov Y)^\Gamma \ra H^2(\Gamma,\ov \Q^*)$$
We have $H^1(\Gamma,\ov \Q^*)=Hom(\Gamma,\ov \Q^*)=\Z/2$,
and $H^2(\Gamma,\ov \Q^*)=0$ by periodicity.
On the other hand, $Ker(\psi^*)=\Z/2$.
The snake lemma applied to the upper part of diagram (2.3)
now gives rise to the exact sequence of $G$-modules
$$0\ra\Z/2\ra\Z/2\ra K \ra J[2]\ra 0.$$
This proves Lemma 1.\  QED
\bigskip
{\bf Corollary.} {\it We have} $\Br(\ov X)^G=0$, {\it thus}
$\Br(X)=\Br_1(X)$.
\medskip
{\it Proof.} This follows from Lemma 1 since 
$\Br(\ov X)^G=Hom_G(J[2](\ov k),\Z/2)=0$.\  QED

\bigskip
{\bf Lemma 2.} {\it The group 
$f^*(\Br(X))/\Br_0(Y)\subset \Br_1(Y)/\Br_0(Y)=
H^1(G,Pic(\ov Y))$ is 
contained in the image of $H^1(\Q,J)[2]\times H^1(\Q,D)$
under the map induced by the natural inclusion
$J(\ov k)\times D(\ov k)=
Pic^0(\ov Y)\hookrightarrow Pic(\ov Y)$.}
\medskip
{\it Proof.} 
By Corollary we have a canonical functorial isomorphism
$${\rm Br}(X)/{\rm Br}_0(X)
\buildrel{=}\over{\hbox to 8mm{\rightarrowfill}} 
%\buildrel{{\tilde}\over{\hbox to 8mm{\rightarrowfill}}} 
H^1(G,Pic(\overline X)).$$  
By functoriality of spectral sequence (2.2) 
we have to consider the subgroup of ${\rm Br}_1(Y)/{\rm Br}_0(Y)$ isomorphic to the 
image of the map 
$$f^*: H^1(G,Pic(\ov X))\ra H^1(G,Pic(\ov Y)).$$ 
Observe that 
$NS(\ov Y)^{\Gamma}\subset NS(\ov Y)$ is torsion 
free. We have $f_*f^*(x)=2x$, for any $x\in NS(\ov X)$, and 
$y+(\sigma, \rho)y=f^*f_*(y)$, for any $y\in NS(\ov Y)$.
Therefore the $G$-modules
$NS(\ov X)\otimes \Q$ and 
$NS(\ov Y)^{\Gamma}\otimes \Q$ are isomorphic. Since $X_\C$ has second Betti number $b_2=2$, we have $dim(NS(\ov Y)^{\Gamma}\otimes \Q)=2$.
The classes of fibres of canonical projections 
$\pi_1: \ov Y\ra \ov C$ and $\pi_2:\ov Y\ra \ov D$ give two linearly
independent $G$-invariant elements of this vector space,
implying that it carries trivial $G$-action. 
Thus $NS(\ov Y)^{\Gamma}=\Z\oplus\Z$ as a $G$-module. 
Consider again diagram (2.3).
Note that on replacing $NS(\ov Y)^\Gamma$ by the image of
$Pic(\ov Y)^\Gamma \ra NS(\ov Y)^\Gamma$ we obtain from (2.3)
a commutative diagram whose middle row is right exact.
This image is a submodule of $NS(\ov Y)^{\Gamma}=\Z\oplus\Z$,
and hence is a free abelian group with trivial $G$-action.
Since
$H^1(G,\Z)=0$ the modified diagram (2.3) gives rise to the 
following commutative diagram with exact rows:
$$\diagram{
H^1(\Q,D') & \ra & H^1(G,Pic(\ov X))&\ra 
& H^1(G,NS(\ov X))\cr
\vfl{}{}{3mm}& & \vfl{}{}{3mm}& & \vfl{}{}{3mm}\cr
H^1(\Q,J[2])\times H^1(\Q,D) & \ra &
H^1(G,Pic(\ov Y)^{\Gamma})&\ra& 0\cr
\vfl{}{}{3mm}& & \vfl{}{}{3mm}& & \vfl{}{}{3mm}\cr
H^1(\Q,J)\times H^1(\Q,D) & \ra & 
H^1(G,Pic(\overline C\times \overline D))&\ra&
H^1(G,NS(\ov Y))\cr
}$$
The statement of Lemma 2 follows from the commutativity of this diagram.\  QED
\medskip
Now we can finish the proof of the theorem.
\medskip
We claim that the group 
$f'^*(\Br(X))/\Br_0(Y')\subset \Br_1(Y')/\Br_0(Y)=
H^1(G,Pic(\ov {Y'}))$ is 
contained in the image of $H^1(\Q,D)$
under the map induced by the natural inclusion
$J(\ov k)\times D(\ov k)=Pic^0(\ov {Y'})
\hookrightarrow Pic(\ov {Y'})$.
After Lemma 2 and by functoriality of 
spectral sequence (2.2) in order to prove this
we only have to remark that
the inverse image map $Pic^0(\ov Y)\ra Pic^0(\ov {Y'})$ 
is multiplication by $(2,1)$ on $J(\ov k)\times D(\ov k)$.
The first statement of (b) now follows.
The theorem is proved. \ QED

\bigskip

J.-L. Colliot-Th\'el\`ene conjectured 
that the Manin obstruction to the Hasse principle for 
zero-cycles of degree one is the only obstruction for 
all varieties over a number field $k$ ([CT], Conj. 1.5 (a),
this statement was also formulated and discussed by 
S. Saito in [Sa], Sect. 8). This would imply 
the existence of a $K$-point on $X$ for some odd 
degree extension $K/\Q$.

\bigskip

{\bf 3. Refinement of the Manin obstruction}

\bigskip
Let us now combine the theory of descent and Manin's treatment of
the Brauer--Grothendieck group to define an ``iterated 
Manin obstruction".
The descent with respect to algebraic tori and its
relation to the algebraic part of the Manin obstruction 
was studied in detail by
J.-L. Colliot-Th\'el\`ene and J.-J. Sansuc [CS]. We 
follow this 
work in a more general context of groups of multiplicative
type, which also covers the classical descent on abelian
varieties.
\medskip
Let $X$ be a variety over 
a number field $k$. A reasonable class of varieties for
which the descent method works well are those
satisfying the condition
$$H^0(\ov X,\G_m)=\ov k^*.\leqno{(3.1)}$$
This condition is satisfied by all proper varieties.

Let $M$ be a finitely generated $G$-module, and let
$S=Hom_{k-groups}(M,\G_m)$ be the $k$-group
of multiplcative type dual to $M$. There are two important cases:
if $M$ is finite, then $S$ is finite, and if $M$ is torsion-free, then
$S$ is an algebraic $k$-torus.

The obstruction that we are going to define is
attached to a $G$-homomorphism
$$\lambda: M\lra Pic(\ov X).$$ 
Let 
$\Br_\lambda(X):=r^{-1}\lambda_*(H^1(G,M))\subset
\Br_1(X)$ (the map $r$ comes from
spectral sequence (2.2)), and define 
$X(\A_k)^{\Br_\lambda}\subset X(\A_k)$ as
the set of adelic points orthogonal to $\Br_\lambda(X)$ with
respect to the Brauer--Manin pairing. 

If $Y$ is an $X$-torsor under $S$, and $\sigma\in H^1(X,S)$,
then the ``twist" $f_\sigma: Y_\sigma\ra X$ is defined
as the $X$-torsor under $S$ whose class is 
$[Y_\sigma]=[Y]-\sigma\in H^1(X,S)$.

For $p:X\ra Spec(k)$ satisfying (3.1) there is the 
following fundamental 
exact sequence ([CS], Thm. 1.5.1):
$$0\ra H^1(k,S)\ra H^1(X,S)\ra Hom_G(M,Pic(\ov X))
\buildrel{\partial}\over{\hbox to 6mm{\rightarrowfill}} 
H^2(k,S)\ra H^2(X,S)\leqno{(3.2)}$$
This is the sequence of low degree terms of the 
spectral sequence
$$Ext^p_k(M,R^qp_*{\bf G}_m)\Rightarrow
Ext^{p+q}_X(p^*M,{\bf G}_m).$$
If $Y$ is an $X$-torsor under $S$, then the image of
the class of $Y$ under 
the map $H^1(X,S)\ra Hom_G(M,Pic(\ov X))$ in (3.2)
is called {\it the type} of $Y$. By (3.2) a torsor of a 
given type is unique up to twist by an element of $H^1(k,S)$.

\medskip
{\bf Proposition.} (a) {\it 
Let $X$ be a smooth geometrically integral variety over 
a number field $k$ satisfying condition (3.1).
Suppose that
$X(\A_k)^{\Br_\lambda}\not=\emptyset$.
Then there exists 
an $X$-torsor $Y$ under $S$ of type $\lambda$, 
$f:Y\ra X$, such that }
$$X(\A_k)^{\Br_\lambda}=\cup_{\sigma\in H^1(k,S)}f_\sigma(Y_\sigma(\A_k)).$$ 
\vskip 2mm

(b) {\it When $X$ is proper there exists a finite subset
$\Sigma\subset H^1(k,S)$ such that}
 $$X(\A_k)^{\Br_\lambda}=\cup_{\sigma\in \Sigma}f_\sigma(Y_\sigma(\A_k)),\ \
X(k)=\cup_{\sigma\in \Sigma}f_\sigma(Y_\sigma( k)).$$ 

\vskip 2mm
(c) {\it Suppose there exists an $X$-torsor $Y$ under $S$ of 
type $\lambda$ such that $Y_\sigma(\A_k) \not=\emptyset$,
then $X(\A_k)^{\Br_\lambda}\not=\emptyset$.}

\bigskip
The proof will be given in the next section. 
\medskip
We define an iterated version of the Manin obstruction
related to $\lambda:M\ra Pic(\ov X)$ as follows. 
Define
$$X(\A_k)^{\lambda}:=\cup_{\sigma\in H^1(k,S)}
f_\sigma(Y_\sigma(\A_k)^\Br).$$
For each $\sigma\in H^1(k,S)$ we have the following 
obvious inclusions:
$$Y_\sigma(k)\subset Y_\sigma(\A_k)^{\Br}
\subset Y_\sigma(\A_k)^{f_\sigma^*(\Br(X))}
\subset Y_\sigma(\A_k).$$
Applying $f_\sigma$ and taking the union for all
$\sigma\in H^1(k,S)$ we obtain
$$X(k)\subset X(\A_k)^{\lambda}
\subset X(\A_k)^{\Br}
\subset X(\A_k)^{\Br_\lambda}\subset X(\A_k).$$
Thus the emptiness of $X(\A_k)^{\lambda}$
is an obstruction to the Hasse principle which is 
a refinement of the Manin obstruction related to the whole
of $\Br(X)$. 

In our example consider $D\in Div(X)$ 
such that the divisor of the function
$g(t)$ on $X$ is $2D$ (cf. (2.1)). Let 
$\lambda:\Z/2\hookrightarrow Pic(\ov X)$ be the inclusion
of the cyclic subgroup generated by $D$.
Then $f:Y\ra X$ is ``a torsor of type $\lambda$" 
(by the local description of torsors, [CS], 2.3, 2.4.1). By
the proof of the assertion (a) of the theorem, we
can take $\Sigma=\{1\}\subset H^1(\Q,\Z/2)=\Q^*/\Q^{*2}$. 
We have $Y(\A_\Q)^\Br=\emptyset$ since
the Manin obstruction is the only obstruction to the 
Hasse principle for curves of genus one. This implies $X(\A_\Q)^\lambda=\emptyset$,  
whereas $X(\A_\Q)^\Br\not=\emptyset$ by assertion (b)
of the theorem.

The proposition implies that if the descent
varieties $Y_\sigma$'s have the property that the Manin obstruction
associated to some finite subgroup of 
${\rm Br}(Y_\sigma)/{\rm Br}_0(Y_\sigma)$
is the only obstruction to the Hasse principle, then there still is
a finite process to decide whether $X(k)$ is empty or not.
This is what happens for varieties $X$ given by (2.1): 
testing the Manin obstruction on $Y$ is
reduced to testing it on $C$ and $D$. By the
classical result of Manin ([M], Thm. 6) 
one is led to consider the conjecturally finite
group ${\cyr X}(J)$, and check whether there is 
$\beta\in{\cyr X}(J)$ such that
the Cassels--Tate pairing $<[C],\beta>$ is not zero.
By the conjectural non-degeneracy of the 
Cassels--Tate pairing ([Mi], Thm. 6.13 (a)) such a $\beta$
exists if and only if $[C]=0$, that is, 
precisely when $C$ has rational points.
If $[C]\not=0$, then 
the Manin obstruction on $C$ is not empty, 
{\it a fortiori\/} the same is true for $Y$.

\bigskip

{\bf 4. Descent with respect to groups of multiplicative 
type} 

\bigskip
In this section we summarize the analysis of relations 
between $X$-torsors under multiplicative groups and
the algebraic part of $\Br(X)$ undertaken in [CS]. Our
treatment is somewhat more general since we do not assume
that the elementary obstruction for the existence of
$k$-points on $X$ vanishes.

The statement {\it (ii)\/}
of the following lemma is a complement to [CS]
clarifying the relation between the fundamental exact
sequence (3.2) 
and the sequence of low degree terms of (2.2). It
also provides an explicit construction of Azumaya algebras
as cup-products via the pairing
$$H^1(k,M)\times H^1(X,S)
\buildrel{(p^*,id)}\over{\hbox to 10mm{\rightarrowfill}} 
H^1(X,M)\times H^1(X,S)\ra\Br(X).$$
Recall that $r:\Br_1(X)\ra H^1(G,Pic(\ov X))$
is the canonical map from spectral sequence (2.2).
\medskip

{\bf Lemma 3.} {\it Let $X$ be a 
variety over a field $k$ satisfying (3.1),
$M$ be a finitely generated $G$-module, $S$ be the dual
$k$-group of multiplicative type. 
Let $\lambda\in
Hom_k(M,Pic(\ov X))$. Suppose there exists an $X$-torsor
$T$ under $S$ of type $\lambda$. Then for $\alpha\in 
H^1(k,M)$ we have
\vskip 2mm
(i) $p^*(\alpha)\cup [T]\in\Br_1(X)$, 
\vskip 2mm

(ii) $r(p^*(\alpha)\cup [T])=\lambda_*(\alpha)$,
\vskip 2mm

(iii) let $A=p^*(\alpha)\cup [T]\in H^2(X,{\bf G}_m)$,
then for $P\in X(k)$ the specialization
$A(P)\in \Br(k)$ equals $\alpha \cup T(P)$ with
respect to the natural pairing $H^1(k,M)\times H^1(k,S)\ra
\Br(k)$,
\vskip 2mm

(iv) for any $A\in\Br_\lambda(X)$ and any $\alpha \in H^1(k,M)$ such that 
$r(A)=\lambda_*(\alpha)$, there exists 
$a_0\in \Br(k)$
such that for any point $P\in X(K)$
over an extension $k\subset K$ we have 
$$Res_{k,K}(\alpha)\cup T(P)=A(P)+Res_{k,K}(a_0)\ 
\ \in\Br(K).$$
}
\medskip
{\it Proof. (i)\/} is trivial. 

{\it (ii)\/} We have 
$Ext^n_X(p^*M,{\bf G}_m)=H^n(X,S)$ ([CS], Prop. 1.4.1), and
also $Ext^n_X(\Z,p^*M)=H^n(X,p^*M)$. 
There are Yoneda cup-products
$$H^1(X,p^*M)\times Ext^n_X(p^*M,{\bf G}_m)\ra 
H^{n+1}(X,{\bf G}_m)\leqno{(4.1)}$$
and
$$H^1(k,M)\times Ext^n_k(M,Pic(\ov X))\ra 
H^{n+1}(k,Pic(\ov X)).\leqno{(4.2)}$$
Let 
$$0\ra M\ra N\ra \Z\ra 0\leqno{(4.3)}$$
be an exact sequence of $G$-modules whose class in 
$Ext^1_k(\Z,M)=H^1(G,M)$ is $\alpha$. 
Let $d$ be the connecting homomorphism in the long
exact sequence of $Ext$'s in the first variable
assotiated to (4.3).
By the definition
of Yoneda cup-product (4.1), for
$\xi\in Ext^n_X(p^*M,{\bf G}_m)$ we have 
$p^*(\alpha)\cup\xi=d(\xi)$. 
In case of (4.2) we have for
$\zeta\in Ext^n_k(M,Pic(\ov X))$ the analogous formula $\alpha\cup\zeta=d(\zeta)$.

We claim that there is the following commutative diagram:
$$\diagram{
Ext_X^1(p^*M,{\bf G}_m) & \lra &Hom_k(M,Pic(\ov X))\cr
\vfl{d}{}{3mm}&&\vfl{d}{}{3mm}\cr
\Br_1(X)&\lra&H^1(G,Pic(\ov X))\cr
}\leqno{(4.4)}$$
Here the upper arrow is the map in (3.1) which associates
to a torsor its type, and the
lower arrow is the map $r$ obtained from the spectral 
sequence (2.2). The statement {\it (ii)} is a consequence of
the commutativity of this diagram.

In order to prove this commutativity it will be
convenient to place ourselves in a 
more general context. Let
$$0\ra A\ra B\ra C\ra 0\leqno{(4.5)}$$
be an exact sequence of $G$-modules.
Let ${\cal D}^+$ be the derived
category of bounded below complexes of $G$-modules, and
$$\ldots\ra C^{\bullet}[-1]\ra A^{\bullet}\ra B^{\bullet}\ra C^{\bullet}\ra \ldots$$
be the corresponding distinguished triangle in ${\cal D}^+$.
For any $F\in {\cal D}^+$, $i\in \Z$, the truncation functors
provide distinguished triangles
$$\ldots\ra  \tau_{\leq i}(F)\ra F\ra \tau_{\geq i+1}(F)\ra \ldots$$
We have a commutative diagram
$$\diagram{
Hom(A^{\bullet},F) & \lra &
Hom(A^{\bullet},\tau_{\geq 1}(F))\cr
\vfl{d}{}{3mm}&&\vfl{d}{}{3mm}\cr
Hom(C^{\bullet},F)[1] & \lra &
Hom(C^{\bullet},\tau_{\geq 1}(F))[1]\cr
}\leqno{(4.6)}$$
Applying this to $F=\tau_{\leq 1}(Rp_*{\bf G}_m)$ and
the sequence (4.3) as (4.5) we obtain the commutative diagram
$$\diagram{
Hom(M^{\bullet},\tau_{\leq 1}(Rp_*{\bf G}_m)) & \lra &
Hom(M^{\bullet},\tau_{\geq 1}(\tau_{\leq 1}(Rp_*{\bf G}_m)))\cr
\vfl{d}{}{3mm}&&\vfl{d}{}{3mm}\cr
\tau_{\leq 1}(Rp_*{\bf G}_m)[1] & \lra &
\tau_{\geq 1}(\tau_{\leq 1}(Rp_*{\bf G}_m))[1]\cr
}$$

On taking the first
hypercohomology groups we get the desired 
diagram (4.4). Indeed, quite generally 
${\bf H}^n(Hom(M^{\bullet},Rp_*{\bf G}_m))=
Ext_X^n(p^*M,{\bf G}_m)$.
This implies that
$${\bf H}^1(Hom(M^{\bullet},\tau_{\leq 1}(Rp_*{\bf G}_m)))=
Ext_X^1(p^*M,{\bf G}_m)$$ since obviously
${\bf H}^0(Hom(M^{\bullet},\tau_{\geq 2}(Rp_*{\bf G}_m)))=
{\bf H}^1(Hom(M^{\bullet},\tau_{\geq 2}(Rp_*{\bf G}_m)))=0$.
Similarly, we have
${\bf H}^1(\tau_{\geq 2}(Rp_*{\bf G}_m))=0,$
and also
${\bf H}^2(\tau_{\geq 2}(Rp_*{\bf G}_m))=
H^2(\ov X,{\bf G}_m)^G,$
hence we get an exact sequence
$$0\ra {\bf H}^2(\tau_{\leq 1}(Rp_*{\bf G}_m))\ra 
H^2(X,{\bf G}_m)\ra H^2(\ov X,{\bf G}_m)^G$$
which identifies ${\bf H}^1(\tau_{\leq 1}
(Rp_*{\bf G}_m)[1])$ with $\Br_1(X)$.
Next, we have 
$$\tau_{\geq 1}(\tau_{\leq 1}(Rp_*{\bf G}_m))=
Pic(\ov X)^\bullet[-1],$$
hence 
$${\bf H}^1(Hom(M^\bullet,
\tau_{\geq 1}(\tau_{\leq 1}(Rp_*{\bf G}_m))))=
{\bf H}^1(Hom(M,Pic(\ov X))^\bullet[-1])=
Hom_k(M,Pic(\ov X)),$$
and similarly
$${\bf H}^1(\tau_{\geq 1}(\tau_{\leq 1}(Rp_*{\bf G}_m))[1])=
{\bf H}^1(Pic(\ov X)^{\bullet})=
H^1(G,Pic(\ov X)).$$
Comparing (4.6) with (4.4) one sees also that the
arrows in these diagrams are identical. The proof of
{\it (ii)\/} is now complete.

{\it (iii)\/} 
This follows from the functoriality of the cup-product.

{\it (iv)\/} Take any $\alpha \in H^1(G,M)$ such that 
$r(A)=\lambda_*(\alpha)$. Then
by {\it (ii)\/} we have that
$$A-p^*(\alpha)\cup [T]\in \Br_0(X),$$ thus it is the image
of some element $a_0\in \Br(k)$. Now apply {\it (iii)\/}.
(Cf. [CS], Lemme 3.5.2.) \ QED

\bigskip
{\it  Proof of the proposition.} (a) and (b)
\medskip

There are three steps:
\medskip
(1) If there exists an adelic point which is 
Brauer--Manin orthogonal
to $\lambda_*({\cyr X}^1(M))\subset \Br_1(X)/\Br_0(X)$, 
then there exists a torsor
$f:Y\ra X$ of type $\lambda$.

(2) If an adelic point is Brauer--Manin orthogonal
to $\Br_\lambda (X)$, then there exists 
$\sigma\in H^1(k,S)$ such that this point lifts to an
adelic point on $Y_\sigma$; the similar assertion for
$k$-points is evident.

(3) If $X$ is proper, then a finite set $\Sigma\subset H^1(k,S)$ would
suffice to cover all points in $X(\A_k)^{\Br_\lambda}$.
\medskip
The statement (1) is very similar to Prop. 3.3.2 of [CS]. 
The case under consideration in {\it loc. cit.} is that
of $Pic(\ov X)$ of finite type and $\lambda: M\ra Pic(\ov X)$
the identity map. We establish (1) along the lines of
that proof. 

It follows from (3.2) that
there exists an $X$-torsor under $S$ of type $\lambda$
if and only if the image $\partial(\lambda)$ of $\lambda$
in $H^2(k,S)$ is zero. 

Let $b_X\in Ext^2_G(Pic(\ov X),\ov k^*)$
be the class of the 2-extension 
$$1\ra \ov k^*\ra \ov k(X)^*\ra Div(\ov X)\ra Pic(\ov X)
\ra 0\leqno{(4.7)}$$
Since $X(\A_k)\not=\emptyset$
this class goes to 0 under the restrictions from
$k$ to $k_v$ ([CS], Prop. 2.2.4, 2.2.2 (b)), so that
$$b_X\in  {\cyr X}(Ext^2_G(Pic(\ov X),\ov k^*)).$$
Recall that $Ext^i_G(M,\ov k^*)=H^i(k,S)$ ([Mi], I.4.12 (a)).
By ([CS], Prop. 1.5.2 {\it iv\/}) we have  
$$\partial(\lambda)=\lambda^*(b_X)\in Ext^2_G(M,\ov k^*).$$ 

One can reinterpret the proof of Lemma 3.3.3 of [CS] as
follows.  Let $A\in \Br_1(X)$ be such that 
$a=r(A)\in {\cyr X}^1(G,Pic(\ov X))$. Therefore $A$
is ``locally constant":  
$i(A):=\sum_{v\in \Omega}\inv_v(A(P_v))$
does not depend on the choice of $\{P_v\}$. By the global
reciprocity $i(A)$ depends only on $a=r(A)$, so that we
adopt the notation $i(a)=i(A)$. Let $a=\lambda_*(\alpha)$
for some $\alpha\in {\cyr X}^1(G,M)$. Then we have the
following formula:
$$i(\lambda_*(\alpha))=<\lambda^*(b_X),\alpha>,
\leqno{(4.8)}$$
where $<,>$ is the canonical duality pairing 
([Mi], I.4.20 (a))
$${\cyr X}^2(G,S)\times{\cyr X}^1(G,M) \ra\Q/\Z$$
To see this one can define a pairing 
${\cyr X}(Ext^2_G(Pic(\ov X),\ov k^*))\times
{\cyr X}^1(G,Pic(\ov X)) \ra\Q/\Z$
following the standard pattern ([Mi], p.79). The
proof of Lemma 3.3.3 of [CS] is precisely checking
that $i(a)$ is given by coupling $b_X$ with $a$ 
with respect to this last pairing. By the functoriality
we get (4.8).

By the assumption of (1) there exists an adelic point 
$\{P_v\}$ such that for any $\alpha\in {\cyr X}^1(G,M)$
we have
$$<\partial(\lambda),\alpha>=
<\lambda^*(b_X),\alpha>=i(\lambda_*(\alpha))=0.$$ 
By the non-degeneracy of the pairing $<,>$ we conclude 
that $\partial(\lambda)=0$. This proves (1).

\medskip
{\it Remark.\/} Note that (4.8) is completely 
analogous to a
result of Manin ([M], Thm. 6). Assume that $X$ is a $k$-torsor under an
abelian variety $S$, $S^t$ is the dual of $S$, $\lambda: 
S^t=Pic^0(\ov X)\hookrightarrow Pic(\ov X)$ is the
natural inclusion. Then presumably 
$\lambda^*(b_X)\in Ext^2_k(S^t,{\bf G}_m)$ coincides
up to a sign with the class 
$[X]\in H^1(k,S)=Ext^2_k(S^t,{\bf G}_m)$, where the last
isomorphism is in ([Mi], I.3.1). If $X(\A_k)\not=\emptyset$,
and $a=r(A)=\lambda_*(\alpha)$ for some 
$\alpha\in {\cyr X}^1(k,S^t)$, then
(4.8) with $<,>$ understood as the canonical
Cassels--Tate
pairing ${\cyr X}^1(k,S)\times {\cyr X}^1(k,S^t)\ra\Q/\Z$ becomes Manin's formula
$\sum_{v\in \Omega}{\rm inv}_v(A(P_v))=<[X],\alpha>$. 

\medskip

Inverting a finite number of primes in the ring of integers 
$\O_k\subset k$ one obtains a ring $\O\subset k$
such that $M$ is unramified over $\O$.
Let $\S=Hom_{\O-groups}(M,\G_m)$.
Let $\O_v$ be the completion of $\O$ at a prime
$v\in Spec(\O)$, $\S_v:={\cal S}\times_{\O}\O_v$.
To siplify the notation we shall write $H^i(R,\cdot)$
for $H^i(Spec(R),\cdot)$ when $R$ is a ring.
We denote $N^*=Hom(N,\Q/\Z)$. The cup-product
$$H^1(k_v,S)\times H^1(k_v,M)\ra \Br(k_v)
\buildrel{{\rm inv}_v}\over{\hbox to 8mm{\rightarrowfill}}
\Q/\Z$$ 
defines an isomorphism (local duality) of finite groups 
$H^1(k_v,S)=H^1(k_v,M)^*$ ([Mi], I.2.3).

Consider the exact sequence (cf. [Mi], I.4.20 (b))
$$0\ra {\cyr X}^1(k,S)\ra H^1(k,S)\ra P^1(k,S)\ra H^1(G,M)^*\leqno{(4.9)}$$
where $P^1(k,S)$ is the restricted product of 
$H^1(k_v,S)$ for all places $v$ with respect to
subgroups $H^1(\O_v,\S_v)$ defined for $v\in Spec(\O)$,
that is, for almost all $v$.
The right arrow in (4.9) is the sum of 
$$\tau_v: H^1(k_v,S)=H^1(k_v,M)^*\ra H^1(k,M)^*$$
which is the composition of the local duality isomorphism and
the dual of the restriction map. Any element
of $H^1(k,M)$ for almost all places $v$ restricts to
an element of $H^1(\O_v,M)\subset H^1(k_v,M)$ ([Mi], I.4.8),
the group orthogonal to $H^1(\O_v,\S_v)$ with respect to
the local duality pairing. Thus the
second arrow in (4.9) is given by a finite sum and hence
is well defined. 

Let $Y$ be an $X$-torsor under $S$ of type $\lambda$ whose
existence was established in (1). Let $Y(P)$ denote
the element in $H^1(k(P),Y)$ given by the fibre of 
$Y$ at a closed point $P$.
The map $X(\A_k)\ra \prod_{v\in \Omega}H^1(k_v,S)$
sending $\{P_v\}$ to $\{Y(P_v)\}$ has its image
in $P^1(k,S)$. (Indeed, up to inverting a finite number of
primes in $\O$ depending on $X$, $Y$, and $\{P_v\}$, 
there exists an
integral, regular model ${\cal X}/{\cal O}$ and
an ${\cal X}$-torsor ${\cal Y}$ under ${\cal S}$
which give $X$, $S$, and $Y$ at the generic point
$Spec(k)$ of $Spec({\cal O})$. Moreover, we can assume that
$P_v\in \X(\O_v)$ for $v\in Spec(\O)$, thus
$\Y(P_v)\in H^1(\O_v,\S_v)$ for such $v$.) 
Thus we can define the map (cf. [CS], Def. 3.4.1)
$$X(\A_k)\ra H^1(k,M)^*, \ {\rm by} \ 
\{P_v\}\mapsto \sum_{v\in \Omega}\tau_v(Y(P_v)).$$
From Lemma 3 {\it (iv)\/} it follows that 
for any $A\in\Br_\lambda(X)$ and any $\alpha\in H^1(G,M)$
such that $r(A)=\lambda_*(\alpha)$ there exists 
$a_0\in \Br(k)$
such that for any $P_v\in X(k_v)$ we have 
$$Res_{k,k_v}(\alpha)\cup Y(P_v)=A(P_v)+Res_{k,k_v}(a_0)\in \Br(k_v).$$
This equality combined with global reciprocity implies
(cf. [CS], Cor. 3.5.3) that
$$\sum_{v\in \Omega} \tau_v(Y(P_v))(\alpha)=\sum_{v\in \Omega}\inv_v(A(P_v)),\leqno{(4.10)}$$
which is just the Brauer--Manin pairing. Let us now complete
the proof of (2). If $\{P_v\}\in X(\A_k)$ is 
Brauer--Manin orthogonal
to $\Br_\lambda (X)$, then we see that
$\{ Y(P_v)\} $ goes to zero under the right arrow
of (4.9) thus is the image of some 
$\sigma\in H^1(k,S)$. Therefore $Y_\sigma(P_v)=0$
for all $v$, and the fibre of $f_\sigma :Y_\sigma \ra X$
over $\{P_v\}$ contains an adelic point, which proves (2).

\medskip

If $X$ is proper one can moreover assume that $\X$ is
proper over $Spec(\O)$. The choice of $\O$, $\X$, and $\Y$
is now independent on $\{P_v\}$. By the properness 
of $X$ we have
$\X(\O_v)=X(k_v)$, and thus
$\Y(P_v)\in H^1(\O_v,\S_v)$ for $v\in Spec(\O)$.
It follows that an element
of $H^1(k,S)$ mapping to $\{Y(P_v)\}\in P^1(k,S)$ 
actually lands in the subset 
$\prod_{v\notin Spec(\O)}H^1(k_v,S)\times
\prod_{v\in Spec(\O)}H^1(O_v,\ S_v)$. 
Let us show that this implies that
$x\in H^1({\cal O},{\cal S})$, the group which 
is well known to be finite. (This is done by reduction
to the case of a finite group and of a torus. The finiteness
in these cases is proved in [M], II.2.13, and II.4.6.)
Let $\F_v$ be the 
residue field at $v$. Consider
the following commutative diagram whose upper row is
the exact sequence of the pair $Spec(k)\subset Spec(\O)$,
the lower row is the one of the pair
$Spec(k_v)\subset Spec(\O_v)$, and the vertical arrows
are the natural restrictions:
$$\diagram{
H^1(\O,{\cal S}) & \lra &H^1(k,S)& \lra &
\prod_{v\in Spec(\O)}H^2_{Spec(\F_v)}(\O,\S)     \cr
\vfl{}{}{3mm}&&\vfl{}{}{3mm}&&\vfl{}{}{3mm}\cr
\prod_{v\in Spec(\O)}H^1(\O_v,\S_v) & \lra &
\prod_{v\in Spec(\O)}H^1(k_v,S)& \lra &
\prod_{v\in Spec(\O)}H^2_{Spec(\F_v)}(\O_v,\S_v)     \cr
}$$
By the \'etale excision theorem combined with Greenberg's
theorem one concludes that the right vertical arrow
is an isomorphism, which proves what we need. (I am
grateful to J.-L. Colliot-Th\'el\`ene for this argument.)
This completes the proof of statement (3).
Hence the parts (a) and (b) of the proposition are
proved. 
\medskip
{\it  Proof of the proposition.} (c)
\medskip
Let $\{P_v\}\in X(\A_k)$ be the image of an adelic point
on $Y$. Then $Y(P_v)=0$, and (4.10) implies that 
$\{P_v\}$ is contained in $X(\A_k)^{\Br_\lambda}$.
\ QED
\bigskip

{\bf Acknowledgement.} The author thanks 
J.-L. Colliot-Th\'el\`ene for 
pointing out a gap in the first draft of the example, and
D. Zagier for bringing the paper [GPZ] to his attention.
I thank S. Siksek for the appendix. This work was done
while the author enjoyed the warm hospitality of
Max-Planck-Institut f\"ur Mathematik in Bonn.

\bigskip
\centerline{\bf Appendix, by S. Siksek}
\medskip
\centerline{\bf  4-descent}
\bigskip

Let us fix our notation. Let $E$ be an elliptic curve 
over a field $k$ of characteristic zero. Let $C$ be a principal homogeneous space
under $E$ over $k$ with the action of $E$ given by
the map $\mu: E\times C\ra C$. Then $C$ is called an
$n$-covering if there is a map $\psi:C\ra E$ such that 
the following diagram commutes (cf. [S], Sect. 2):
$$\diagram{
E\times C &\hfl{(n,\psi)}{}{10mm}&E\times E\cr
\vfl{\mu}{}{4mm}&&\vfl{+}{}{4mm}\cr
C&\hfl{\psi}{}{10mm}&E\cr
}$$
An $nm$-covering $\psi':C'\ra E$ is a lifting of $\psi:
C\ra E$ if the analogous diagram for $C'$ factorizes
through this one. 
A principal homogeneous space
$C$ under $E$ can be endowed  with a structure of $n$-covering if and only if its class $[C]\in H^1(k,E)$
belongs to the $n$-torsion. This is precisely when 
$[C]$ can be lifted to a cocycle of $H^1(k,E[n])$ (see [C], III, Sect. 2).
By ([C], III, Sect. 4) the $nm$-covering 
$\psi':C'\ra E$ is a lifting of the $n$-covering
$\psi:C\ra E$ if and only if the corresponding
cocycles are related by the map $m_*:H^1(k,E[mn])\ra
H^1(k,E[n])$. This implies that $m[C']=[C]\in H^1(k,E)$.

We shall consider more closely the case
when $n=m=2$ and $C$ is given by
$$y^2=ax^4+cx^2+dx+e.$$ 
Let $K$ be the $k$-algebra $k[x]/(ax^4+cx^2+dx+e )$, and 
let $\theta$ be the image of $x$ in $K$. Let 
$C'\in\P^3_k$ be the intersection
of two quadrics obtained by equating the coefficients
of $\theta^2$ and $\theta^3$ in the following
formula
$$ X-\theta Z = \epsilon (x_1+x_2 \theta+x_3 \theta^2+x_4 \theta^3)^2, \leqno{(A.1)}$$
for some $\epsilon\in K^*$ such that 
$a^{-1}N_{K/k}(\epsilon)\in k^{*2}$. The following
lemma is well known and is given here for the sake of
completeness, cf. [B/SwD], [MSS], [Ca].
\medskip
{\bf Lemma.} {\it (a) $C$ can be equipped with a
structure of 2-covering $\psi$ of
the elliptic curve $E: u^2=v^3-27Iv-27J$, where
$I=12ae+c^2$, $J=72ace-27ad^2-2c^3$, such that $\psi^{-1}(0)\subset C$ is given by $y=0$.

(b) $C'$ can be equipped with a structure of 4-covering 
of $E$ which is a lifting of $\psi:C\ra E$.

(c) If $C(k)=\emptyset$, then $[C']$ is of exact order 4.
}
\medskip
Note
that when $a=1$ the passage from the equation of $C$ to
that of $E$ is precisely the passage from a quartic
polynomial to its cubic resolvent. 

Recall the following well known geometric observation.
Let $X\subset \P^4_k$ be the intersection of two 
quadrics given by $Q(\x)=Q'(\x)=0$. The corresponding
discriminantal curve $Y$ defined by $\mu^2=det(\lambda Q+Q')$
is naturally identified with $Pic^2(X)$. Indeed, 
the curve $Y$
parametrizes families of lines on the quadrics of the pencil
spanned by $Q$ and $Q'$. Then $Y$ is identified with $Pic^2(X)$ by the map
which associates to a family of lines on a quadric the 
divisor class of the intersection of any line of 
this family with $X$. The curve $Y$ comes with the map
$\xi:X\ra Y=Pic^2(X)$ sending a point $P$ to the divisor 
class of $2P$ (geometrically this is the family of lines
containing the tangent to $X$ at $P$).
\medskip

{\it Proof of the lemma.\/} 
(a) The curve $C$ is isomorphic to
the following intersection of two quadrics in $\P_k^4$:
$$Q=ut-x^2, \ \ Q'=-y^2+au^2+cut+dxt+et^2.$$
One computes that 
$$det(\lambda Q-Q')=(1/4)(\lambda^3-2c \lambda^2+(c^2-4ae) \lambda+
ad^2) ,$$
which gives the equation of $E$ after
an obvious change of variables. Thus $Y=E$. Let us choose 
the origin of
the group law of $E$ at $\lambda =\infty$ (this point
corresponds to the hyperelliptic divisor on $C$). The map  
$\psi=\xi$ is a 2-covering, and it is easy to see that $\psi^{-1}(0)$ is as required.

(b) Let $\theta_i\in \ov k$ be the roots of the 
equation $ax^4+cx^2+dx+e=0$. 
Let $\epsilon_i\in \ov k$ be the image of $\epsilon$
under the map $K\ra \ov k$ which sends $\theta$ to
$\theta_i$.
Let us denote by
$\delta_i$ the Van der Monde determinant associated to
$\theta_1,\ldots,\theta_4$ with the exception of $\theta_i$,
multiplied by $(-1)^i$. Let 
$$\Delta=
\prod_{1\leq i<j\leq 4}(\theta_i-\theta_j)^2\ \in k$$ 
be the discriminant of  $a^{-1}(ax^4+cx^2+dx+e)$. 

In the pencil of quadrics containing $C'$ we choose 
the following two quadrics individually 
defined over $k$:
$$Q(\x)=\sum_{i=1,2,3,4}\delta_i\epsilon_i
(\sum_{j=1,2,3,4}x_j\theta_i^j)^2,$$
$$Q'(\x)=\sum_{i=1,2,3,4}\theta_i\delta_i\epsilon_i
(\sum_{j=1,2,3,4}x_j\theta_i^j)^2.$$
One computes that 
$$det(\lambda Q-Q')=N_{K/k}(\epsilon)\Delta^2\prod_{i=1,2,3,4}(\lambda-
\theta_i).$$ 
Thus $Y=C$, and we have the map $\xi:C'\ra C=Pic^2(C')$.
A structure of 4-covering $\psi':C'\ra E$ 
is defined by the map
sending a point $P$ to the divisor class of $4P$, and
identifying $Pic^4(C')$ with $E$ by choosing the hyperplane
section class as the origin of the group law. Then $\psi'=\psi\circ\xi$, so that $\psi'$ is a lifting of $\psi$. See ([W], App. III, [MSS], [Ca]) for more on this classical subject. 

(c) Since $\psi'$ is a lifting of $\psi$ 
we have $2[C']=[C]\in H^1(k,E)$. When $C$ has no $k$-point,
$[C]\not=0$, and our statement is proved. \ QED

Let us apply this to the curve
$$ C: \ y^2 = g(x), \ \ g(x)=3(x^4 - 54x^2 - 117x - 243)
\leqno{(A.2)}$$
defined over $\Q$. Using (a) one computes that this is a 2-covering of 
$$J: \ y^2 =x^3 -1221. $$
($C$ is in fact everywhere locally soluble, and was 
initially found using Cremona's program 
{\tt mwrank} [Cr]. By  
([B/SwD], Lemmas 1 and 2, [C], IV, Thm. 1.3) 
any 2-covering which is everywhere locally soluble
can be given by a double cover of ${\bf P}^1_{\bf Q}$.)
It is computed in [GPZ] that $J$ has analytic rank zero.
The authors point out that
the rank of $J(\Q)$ is unconditionally $0$ by the work of Rubin and Kolyvagin. The classical computation of
torsion of such curves implies 
that $E(\Q)=\left\{ 0 \right\}$.
This together with irreducibility of $g(x)$ implies that
$C(\Q)=\emptyset$. 
It is also computed in [GPZ] that $\Sha(J)$ is
predicted to be isomorphic to $\Z/4\times \Z/4$
by  the conjec\-tures of Birch and Swinnerton-Dyer. We shall now exhibit an 
everywhere locally solvable 4-covering 
$C'\ra J$ which is a lifting of the 2-covering $C\ra J$. 
By (c) of the above lemma this gives an element of order 
exactly $4$ in $\Sha(J)$.

It is observed that the element
$\epsilon=- \theta^3/3 - \theta^2 + 29 \theta + 27\in K=
\Q(\theta)$, where $g(\theta)=0$,
has norm $243=3 \times 9^2$. 
We construct a 4-covering $C'$ as
the intersection of two quadrics (A.1), which can be 
written down explicitly:
$${\bf x} A {\bf x}^t = 0 ,~~~{\bf x} B {\bf x}^t=0 $$
where ${\bf x}=(x_1,x_2,x_3,x_4)$ and the entries
of $A$ and $B$ are respectively

\centerline{
$ \diagram{-1 & 11 & -66 & 396 \cr
11 & -66 & 396 & -2520 \cr
-66 & 396 & -2520 & 16335\cr
396 & -2520 & 16335 & -105786
\cr}$\ \ \ 
and\ \ \
$ \diagram{
-1 & -3 & 33 & -198\cr
-3 & 33 & -198 & 1188\cr
33 & -198 & 1188 & -7560\cr
-198 & 1188 & -7560 & 49005
\cr}$}

Let us show that $[C']\in \Sha(J)$.
Criteria for
testing intersections of 2 quadrics in $\P^3$ for everywhere local solubility are given in [MSS]. By
([MSS], Lemma 7) we know that it is soluble  over $\R$, and by ([MSS], Thm. 4) that it is necessary to
test for solubility only at the finite primes $2,~3,~11,~37$
(2 and the divisors of 1221). This result also tells us that the following
will lift to $p$-adic points on $C'$ for $p=2,~3,~11,~37$ respectively:
$(0,2,1,0) mod(2^3),
(12,21,1,0) mod(3^3),
(0,1,0,0) mod(11)
(0,1,9,16) mod(37)$.
This completes our proof that $C'$ represents an element of order exactly $4$
in the the Tate--Shafarevich group of $J$.

We would like to draw the reader's attention to an explicit
form of the Cassels--Tate bilinear pairing on the 2-Selmer
group given in [Ca]. This could have been used to show 
that the Tate--Shafarevich
group does contain an element of exact order $4$, 
though that method would not
have allowed us to write down the equations for
the curve representing such an element.

We thank J. Cremona and N. Smart for their comments on an 
earlier version of this appendix.

\bigskip
\centerline{\bf References}
\medskip
[B] A. Beauville. {\it Surfaces algebriques complexes.} Ast\'erisque
{\bf 54} (1978).

[B/SwD] B. J. Birch and H. P. F. Swinnerton-Dyer. {\it Notes on
elliptic curves. I.} J. reine angew. Math. {\bf 212} (1963) 7-25.

[C] J. W. S. Cassels. {\it Arithmetic on Curves of Genus 1.}
{\it III. The Tate-Shafarevich
and Selmer Groups.\/} Proc. London Math. Soc. {\bf 12} (1962)
 259-296;
{\it IV. Proof of the Hauptvermutung.\/}
J. reine angew. Math. {\bf 211} (1962) 95-112.

[Ca] J. W. S. Cassels. {\it Second descents for elliptic curves.\/} To appear.

[CT] J.-L. Colliot-Th\'el\`ene. {\it L'arithm\'etique 
du groupe de Chow des z\'ero-cycles.} J. th\'eorie 
des nombres de Bordeaux {\bf 7} (1995) 51-73.

[CS] J.-L. Colliot-Th\'el\`ene et J.-J. Sansuc. {\it La
descente sur les vari\'et\'es rationnelles. II.} Duke
Math. J. {\bf 54} (1987) 375-492.

[AA] J.-L. Colliot-Th\'el\`ene, A. N. Skorobogatov and Sir Peter
Swinnerton-Dyer. {\it Double fibres and double covers: paucity of
rational points.} Acta Arithm. {\bf 79} (1997) 113-135.

[Cr] J. E. Cremona. {\it Algorithms for Modular Elliptic Curves.\/}
Cambridge University Press, 1992.

[GPZ] J. Gebel, A. Peth\H{o}, H. G. Zimmer. {\it On Mordell's Equation.} Comp. Math., to appear.

[G] A. Grothendieck. {\it Le groupe de Brauer. I, II, III. \/}
Dans: {\it Dix expos\'es sur la cohomologie des sch\'emas.} 
A. Grothendieck, N. H. Kuipers, Eds. North-Holland (1968)
46-188.

[M] Yu. I. Manin. {\it Le groupe de Brauer-Grothendieck en
g\'eometrie diophantienne.} Actes Congr\`es Int. Math. Nice 1970, tome {\bf I}, Gauthier-Villars, Paris (1971) 401-411.

[MSS] J. R. Merriman, S. Siksek and N. P. Smart, {\it
Explicit $4$-Descents on an Elliptic Curve.}
Acta Arith. {\bf 77} (1996) 385--404.
 
[Mi] J. S. Milne. {\it Arithmetic Duality Theorems.} 
Persp. Math. {\bf 1}, Academic Press, Boston (1986).

[Mi2] J. S. Milne. {\it \'Etale Cohomology.} Princeton Univ.
Press (1980).

[Sa] S. Saito. {\it Some observations on motivic cohomology
of arithmetic schemes.\/} Inv. Math. {\bf 98} (1989) 371-404.

[SW] P. Sarnak and L. Wang. {\it Some hypersurfaces in ${\bf P}^4$
and the Hasse-principle.} C. R. Acad. Sci. Paris {\bf 321} (1995) 319-322.

[Sh] I. R. Shafarevich {\it et al.} {\it Algebraic surfaces.} Proc.
Steklov Inst. Math. {\bf 75} (1967).

[S] Sir Peter Swinnerton-Dyer. {\it Rational points
on certain intersections of two quadrics.\/} In: 
{\it Abelian Varieties.\/} Barth, Hulek, Lange, Eds. Walter
de Gruyter (1995) 273-292.

[W] A. Weil. {\it Number Theory. An approach through 
history.\/} Birkh\"auser, Boston (1983).
\bigskip
\bigskip
Institute for Problems of Information Transmission

Russian Academy of Sciences

19 Bolshoi Karetnyi

Moscow 101447 RUSSIA
\end